\documentclass[a4paper,12pt]{article}

\usepackage{amsmath,amsfonts,amssymb,bm,amsthm}
\usepackage{mathrsfs}
\usepackage{cases}
\usepackage{cite}
\usepackage{indentfirst}
\usepackage[usenames]{color}
\usepackage{hyperref}
\usepackage{fancyhdr}
\usepackage{geometry}
\definecolor{webgreen}{rgb}{0,.5,0}
\definecolor{webbrown}{rgb}{.6,0,0}

\usepackage{color}

\newtheorem{theorem}{Theorem}
\newtheorem{lemma}{Lemma}

\newtheorem{corollary}{Corollary}

\allowdisplaybreaks

\begin{document}

\title{Some extensive discussions of Liouville's theorem and Cauchy's integral theorem on structural holomorphic }

\author{ Gen Wang \thanks{ wanggen@zjnu.edu.cn}}

\date{\em{\small{Department of Mathematics, Zhejiang Normal University,\\ ZheJiang, Jinhua 321004, China.}} }

\maketitle

\begin{abstract}
Classic complex analysis is built on structural function $K=1$ only associated with Cauchy-Riemann equations, subsequently various generalizations of Cauchy-Riemann equations start to break this situation. The goal of this article is to show that only structural function $K=Const$ such that Liouville's theorem is held, otherwise, it's not valid any more on complex domain based on structural holomorphic, the correction should be $w=\Phi {{e}^{-K}}$, where $\Phi =Const$.  Those theories in complex analysis which keep constant are unable to be held as constant in the framework of structural holomorphic.  Synchronously, it deals with the generalization of Cauchy's integral theorem by using the new perspective of structural holomorphic, it is also shown that some of theories in the complex analysis are special cases at $K=Const$,  which are narrow to be applied such as maximum modulus principle.


\end{abstract}

\tableofcontents
\section{Introduction}
As all known, Liouville's theorem and maximum modulus principle all describe the constant properties on different conditions, Cauchy's integral theorem also can be generalized as well\cite{2,5,6,11}. On structural holomorphic \cite{1} of them might be modified to suit general situation.

\subsection{Structural holomorphic}
In the field of complex analysis, the Cauchy-Riemann equations as a certain continuity and differentiability criteria for a complex function $w$ to be complex differentiable, namely, holomorphic or analytic\cite{1,3}.
\begin{equation}
  \left( \begin{matrix}
   \frac{\partial }{\partial x} & -\frac{\partial }{\partial y}  \\
   \frac{\partial }{\partial y} & \frac{\partial }{\partial x}  \\
\end{matrix} \right)\left( \begin{matrix}
   u  \\
   v  \\
\end{matrix} \right)=0\Leftrightarrow \frac{\partial w}{\partial \overline{z}}=0
\end{equation}
Subsequently,  the Carleman-Bers-Vekua (CBV) equation follows\cite{7,8,9,10}
\[\frac{\partial w}{\partial \overline{z}}+Aw+B\overline{w}=0\]where   $A=\frac{1}{4}\left( a+d+\sqrt{-1}c-\sqrt{-1}b \right),~~B=\frac{1}{4}\left( a-d+\sqrt{-1}c+\sqrt{-1}b \right)$,
or equivalently in a matrix form
\begin{equation}\label{eq8}
  \widehat{O}\left( \begin{matrix}
   u  \\
   v  \\
\end{matrix} \right)=0
\end{equation}
where $\widehat{O}=\left( \begin{matrix}
   \frac{\partial }{\partial x}+a & -\left( \frac{\partial }{\partial y}-b \right)  \\
   \frac{\partial }{\partial y}+c & \frac{\partial }{\partial x}+d  \\
\end{matrix} \right)$ is operator matrix.

Furthermore, the boundary value problem was generally investigated:
\[\left\{ \begin{matrix}
   {{w}_{\overline{z}}}+A\left( z \right)w+B\left( z \right)\overline{w}=H,\left| z \right|<1  \\
   \operatorname{Re}\left[ {{z}^{-k}}w\left( z \right) \right]=g\left( z \right),\left| z \right|=1  \\
\end{matrix} \right.\]
What's more, the nonlinear Cauchy-Riemann (NCR) equations is proposed \cite{4}
\[{{u}_{y}}=-{{v}_{x}}+f\left( u,v \right),{{u}_{x}}={{v}_{y}}+g\left( u,v \right)\]
The matrix form is
\[\left( \begin{matrix}
   \frac{\partial }{\partial x} & -\frac{\partial }{\partial y}  \\
   \frac{\partial }{\partial y} & \frac{\partial }{\partial x}  \\
\end{matrix} \right)\left( \begin{matrix}
   u  \\
   v  \\
\end{matrix} \right)=\left( \begin{matrix}
   g\left( u,v \right)  \\
   f\left( u,v \right)  \\
\end{matrix} \right)\]
where $f(u,v)$ and $g(u,v)$ are given functions.
the system of inhomogeneous Cauchy-Riemann equations is
\[{\frac {\partial u}{\partial x}}-{\frac {\partial v}{\partial y}}=G (x,y),~~{\frac {\partial u}{\partial y}}+{\frac {\partial v}{\partial x}}=F (x,y)\]
for some given functions $G(x,y)$ and $F(x,y)$ defined in an open subset of $\mathbb{R}^{2}$.  Structural holomorphic can be expressed in the form \[{{a}_{1}}{{u}_{y}}=-{{a}_{2}}{{v}_{x}}+f\left( u,v \right),~~{{a}_{3}}{{u}_{x}}={{a}_{4}}{{v}_{y}}+g\left( u,v \right)\]
where ${{a}_{i}},i=1,2,3,4$ are real function.
the structural holomorphic can be rewritten \[\frac{\partial w}{\partial \bar{z}}+w\frac{\partial \kappa}{\partial \bar{z}}=0\]
Above all, one has retrospected the developing history of relating holomorphic conditions for a complex function $w$ defined on the complex domain $\mathbb {C}$. In conclusion, the evolution of Cauchy-Riemann equations can be shown in a chain
\begin{align}
  & \frac{\partial w}{\partial \bar{z}}=0\to \frac{\partial w}{\partial \bar{z}}=\varphi \left( z,\overline{z} \right)\to \frac{\partial w}{\partial \bar{z}}+Aw=0 \notag\\
 & \to \frac{\partial w}{\partial \bar{z}}+Aw+B\overline{w}=0\to K\frac{\partial w}{\partial \bar{z}}+w\frac{\partial K}{\partial \bar{z}}=0 \to  \frac{\partial w}{\partial \bar{z}}+w\frac{\partial K}{\partial \bar{z}}=0  \notag
\end{align}Note that it's a trivial term $K\frac{\partial w}{\partial \bar{z}}$ for this $K$ in the $K\frac{\partial w}{\partial \bar{z}}+w\frac{\partial K}{\partial \bar{z}}=0$, then we only study  structural holomorphic  $\frac{\partial w}{\partial \bar{z}}+w\frac{\partial K}{\partial \bar{z}}=0 $.  And more general formally equation can be
\[\frac{\partial w}{\partial \bar{z}}+w\frac{\partial K}{\partial \bar{z}}=\varphi \left( z,\overline{z} \right)\]
where $K$ is structural function which can be chosen arbitrarily, $\varphi \left( z,\overline{z} \right)$ is a given function, it becomes the structural holomorphic as $\varphi \left( z,\overline{z} \right)=0$.

Eventually, most of equations for holomorphic condition are not specified, especially, the coefficients in the equation are indistinct, particularly, the theories for the generalization of Cauchy-Riemann equation most of them come from other branches such as partial differential equation, or physical problems which are not from the pure complex analysis. Surely, the generalizations of Cauchy-Riemann equation also get new difficulties and complexity emerged.  On account of this point, \cite{1} has restarted with the continuity and differentiability of a complex function $w$ as usually done to discuss in a $K$-transformation method given by \[w\left( z \right)\to \widetilde{w}\left( z \right)=w\left( z \right)K\left( z \right)\]
where structural function $K\left( z \right)={{k}_{1}}+\sqrt{-1}{{k}_{2}}$ is complex valued function on $\Omega$, ${{k}_{1}},{{k}_{2}}$ are real functions with respect to the variables $x,y$.  Mathematically, this $K$-transformation is the most natural way to generalize the Cauchy-Riemann equation in complex analysis.

  As traditionally depicted as $K=1$, or the parts take values ${{k}_{1}}=1,{{k}_{2}}=0$.  Mathematically, the generalized holomorphic condition can be clearly expressed by structural function $K\left( z \right)$, and all equations for developing the Cauchy-Riemann equation can be summarized as some of inevitable connections to the structural function $K\left( z \right)$ here. Therefore it follows the main theorem given by
\begin{theorem}\cite{1}\label{t3}
  Let $\Omega\subset \mathbb {C}^{n}$ be an open set and complex valued function $w\left( z \right)=u+\sqrt{-1}v$ is said to be a structural holomorphic on $\Omega$ if and only if
\begin{equation}\label{eq1}
  \frac{{\rm{D}}}{\partial \overline{z^{i}}}w=\frac{\partial w}{\partial \overline{z^{i}}}+w\frac{\partial K}{\partial \overline{z^{i}}}=0 ,~~i=1,2,\cdots ,n
\end{equation}
and its solutions are called structural analytic functions.
\end{theorem}
As pictured in \cite{1}, structural holomorphic differential equation can unify all situations on the way to develop and generalize the Cauchy-Riemann equation in a simple and compact differential form, since $K\left( z \right)$ can be chosen arbitrarily.  Thusly, it suits either functions of one complex variable $z$ in $\mathbb {C}$ or functions of several complex variables $z={{\left( {{z}^{1}},{{z}^{2}},\cdots ,{{z}^{n}} \right)}^{T}}\in {{\mathbb{C}}^{n}}$.
Hence, if one consider \eqref{eq1} in $\Omega\subset \mathbb {C}^{1}=\mathbb {C}$, then \eqref{eq1} is equivalent to the following matrix equation
\[\left( \begin{matrix}
   {{k}_{1}} & {{k}_{2}}  \\
   -{{k}_{2}} & {{k}_{1}}  \\
\end{matrix} \right)\left( \begin{matrix}
   \frac{\partial }{\partial x} & \frac{\partial }{\partial y}  \\
   -\frac{\partial }{\partial y} & \frac{\partial }{\partial x}  \\
\end{matrix} \right)\left( \begin{matrix}
   v  \\
   u  \\
\end{matrix} \right)=\left( \begin{matrix}
   \frac{\partial }{\partial x} & \frac{\partial }{\partial y}  \\
   -\frac{\partial }{\partial y} & \frac{\partial }{\partial x}  \\
\end{matrix} \right)\left( \begin{matrix}
   -{{k}_{2}} & -{{k}_{1}}  \\
   -{{k}_{1}} & {{k}_{2}}  \\
\end{matrix} \right)\left( \begin{matrix}
   u  \\
   v  \\
\end{matrix} \right)\]
Notice that here $\left( \begin{matrix}
   {{k}_{1}} & {{k}_{2}}  \\
   -{{k}_{2}} & {{k}_{1}}  \\
\end{matrix} \right)\left( \begin{matrix}
   \frac{\partial }{\partial x} & \frac{\partial }{\partial y}  \\
   -\frac{\partial }{\partial y} & \frac{\partial }{\partial x}  \\
\end{matrix} \right)$ and $\left( \begin{matrix}
   \frac{\partial }{\partial x} & \frac{\partial }{\partial y}  \\
   -\frac{\partial }{\partial y} & \frac{\partial }{\partial x}  \\
\end{matrix} \right)\left( \begin{matrix}
   -{{k}_{2}} & -{{k}_{1}}  \\
   -{{k}_{1}} & {{k}_{2}}  \\
\end{matrix} \right)$ are structural matrix operators only associated with the Wirtinger derivatives and structural function $K\left( z \right)$. As a result, since $K\left( z \right)$ can be arbitrarily chosen, then ${k}_{1}$ and ${k}_{2}$ are chosen arbitrarily.

To assume that $K\left( z \right)=1+\kappa \left( z \right)$ in $\mathbb {C}$ which is prominently showing the classic complex analysis built on $1$. The structural holomorphic is
\[\frac{\partial w}{\partial \overline{z}}=-w\frac{\partial \kappa }{\partial \overline{z}}\]
Conversely, it laterally reveals the importance of structural function to enlarge the scope of application. \eqref{eq1} can be written in a more general form  \[ \frac{{\rm{D}}}{\partial \overline{z^{i}}}w=\frac{\partial w}{\partial \overline{{{z}^{i}}}}+w\frac{\partial K}{\partial \overline{{{z}^{i}}}}=\varphi \left( z,\overline{z} \right)\]
where $\varphi \left( z,\overline{z} \right)$ is some given function.

Here one gives a classical result below as a reference
\begin{lemma}[\cite{2,6}]\label{l1}
  If $w:\Omega \to \mathbb{C}$ is analytic, then $w$ preserves angles at each point $z_{0}$ of $\Omega$, where $w'\left( {{z}_{0}} \right)\ne 0$.
\end{lemma}
This lemma will tell us that $w$ preserves angles at each point $z_{0}$ of $\Omega$ in the structural holomorphic.

Accordingly, a natural corollary can be deduced in $\mathbb {C}^{n}$,

\begin{corollary}\label{c1}
 The structural holomorphic for  $w\left( z \right)=u+\sqrt{-1}v$ on $\Omega\subset \mathbb {C}^{n}$ is
\[\overline{{\rm{D}}}w=\overline{\partial }w+w\overline{\partial }K=0\]where $K\left( z \right)$ is structural function defined on $\Omega$.
\end{corollary}
If $K=c\in \mathbb{C}$ takes values as complex constant, then $c\overline{\partial }w=0$ always holds, in particular, let $K=1$ be given for $\overline{\partial }w=0$.  In the same way,  $K\left( z \right)=1+\kappa \left( z \right)$ holds for  \[\overline{\partial }w=-w\overline{\partial }\kappa \]
It strongly implies the necessary of the structural function $\kappa \left( z \right)$ for any generalization of Cauchy-Riemann equation.

Let all structural holomorphic function be a collection, denoted as $Shol$.
Note that corollary \ref{c1} is complete and compact form including all possible structural holomorphic condition because of $K\left( z \right)$ can be chosen arbitrarily, say, given a point ${{z}_{0}}\in \Omega \subset {{\mathbb{C}}^{n}}$,  then
$\overline{{\rm{D}}}w\left| _{z={{z}_{0}}}=0 \right.$, that is
\[\overline{\partial }w\left| _{z={{z}_{0}}} \right.+w\left( {{z}_{0}} \right)\overline{\partial }K\left| _{z={{z}_{0}}} \right.=0\]Suppose that $w\left( {{z}_{0}} \right)=0$, then the equation
\[\overline{\partial }w\left| _{z={{z}_{0}}} \right.=0\]holds at ${{z}_{0}}$, accordingly, there is $\overline{\partial }w\left| _{z={{z}_{0}}} \right.=0$ reduced to the classic Cauchy-Riemann equation.

To consider the ${{z}_{0}}\in \Omega \subset {{\mathbb{C}}}$, then \eqref{eq1} becomes $\frac{{\rm{D}}}{\partial \overline{z}}w=\frac{\partial w}{\partial \overline{z}}+w\frac{\partial K}{\partial \overline{z}}=0$ as $n=1$ with  $\frac{{\rm{D}}}{\partial \overline{z}}w=\frac{\partial }{\partial z}w+w\frac{\partial K}{\partial z}$, hence $\frac{\partial }{\partial z}w\neq0$ holds at ${{z}_{0}}$ in $\Omega$, according to the lemma \ref{l1}, structural holomorphic always preserves angles except some particular situations,

A natural question arises while investigating the above-mentioned differential systems structural holomorphic condition: how should the solutions be obtained. In fact, the solution of  \eqref{eq1} in $\mathbb {C}$ or $\mathbb {C}^{n}$ is a little hard to get owing to the specific expression of $K\left( z \right)$ is unknown, but for some particular circumstances, one can obtain the formal solution.

For a reason that $K\left( z \right)$ can be arbitrarily chosen, for instance, taking the form of $K$ is  \[K\left( z \right)={{\overline{z}}^{n}}+{{a}_{1}}{{\overline{z}}^{n-1}}+{{a}_{2}}{{\overline{z}}^{n-2}}+\cdots +{{a}_{n-2}}{{\overline{z}}^{2}}+{{a}_{n-1}}\overline{z}+{{a}_{n}}\]then structural holomorphic condition is
\begin{align}
  & \frac{\partial w}{\partial \overline{z}}=-w\left( n{{\overline{z}}^{n-1}}+{{a}_{1}}\left( n-1 \right){{\overline{z}}^{n-2}}+{{a}_{2}}\left( n-2 \right){{\overline{z}}^{n-3}}+\cdots +2{{a}_{n-2}}\overline{z}+{{a}_{n-1}} \right) \notag
\end{align}
where ${{a}_{i}}\in \mathbb{C}$, or ${{a}_{i}}\in \mathbb{R}$,  $i=1,2,\cdots ,n$, hence choose origin $z=0$ such that structural holomorphic condition turns to  $\frac{\partial w}{\partial \bar{z}}+w{{a}_{n-1}} =0$.

For example, if let
$K\left( z \right)=\overline{z}$ be given in $\mathbb {C}$, it follows \[\frac{\partial w}{\partial \overline{z}}=-w\]or if case
$K\left( z \right)=1+z\sin \overline{z}$ holds, then structural holomorphic condition is given
\[\frac{\partial w}{\partial \overline{z}}+zw\cos \overline{z}=0\]More typical expressions of $K$ can be chosen and considered \[K\left( z \right)=\frac{1}{\sin z},K\left( z \right)=\sin \frac{1}{z},K\left( z \right)=\frac{{{z}^{3}}-2z+10}{{{z}^{5}}+3z-1}\]or \[K\left( z \right)=\frac{{{e}^{z}}}{z},K\left( z \right)=\frac{\sin z}{{{\left( z-1 \right)}^{2}}},K\left( z \right)={{e}^{^{\frac{1}{z}}}},K\left( z \right)=\frac{1}{\sin \left( \frac{1}{z} \right)}\]and more familiar expressions such as $K\left( z \right)={{e}^{z}},\cos z,\sin z,\ln \left( 1+z \right),{{\left( 1+z \right)}^{\alpha }},\frac{1}{1-z}$ on the whole complex plane, all this examples emphasize the arbitrariness of structural function $K\left( z \right)$, then generalized structural Wirtinger derivatives can follow
$ \frac{{\rm{D}}}{\partial \overline{z}}=\frac{\partial }{\partial z}+\frac{\partial K}{\partial z}$.
To one's surprised, one found that structural holomorphic condition remains the same when $K\left( z \right)=\frac{{{e}^{z}}}{z}$ or $\overline{z}$, this phenomenon implies that different structural function $K\left( z \right)$ can derive the same structural holomorphic condition.

On the foundation of above theorem and corollary, one will use them to reanalyze some fundamental theorems within functions of one complex variable $z$ in $\mathbb {C}$, it indicates that most of theorems are not suitable for additional structure in complex domain.  Analytic or holomorphic
is denoted by set $Hol$.

\subsection{ Liouville's theorem }


\begin{lemma}\label{l2}
Suppose that $w$ is a holomorphic function on a connected open set
$\Omega \subset \mathbb{C}$. If $\frac{\partial w}{\partial z}=0$ on $\Omega$, then $w={{e}^{C_{1} }}$ is constant on $\Omega$, where  $C_{1}\in \mathbb{C}$ is a complex constant.
\begin{proof}
  Due to the fact that holomorphic function $w$ satisfies equation $\frac{\partial }{\partial \overline{z}}w=0$ in an open set
$\Omega \subset \mathbb{C}$, together with the extra condition $\frac{\partial w}{\partial z}=0$ on $\Omega$, namely, the systems of linear differential equation $\left\{ \begin{matrix}
   \frac{\partial }{\partial \bar{z}}w=0  \\
   \frac{\partial }{\partial z}w=0  \\
\end{matrix} \right.$,  then it's not difficult to prove that $w={{e}^{C_{1} }}$ is constant on $\Omega$, where $C_{1}\in \mathbb{C}$ is a complex constant.
\end{proof}

\end{lemma}
Hence, with the assistance of the above lemma, if function $w$ is holomorphic defined on the whole complex domain $\Omega$, it's said to be entire function. If plus additional condition for the $w$, then one can discuss more properties about it such as $w$ bounded on $\Omega$, it leads to the well-known Liouville's theorem.  An analytic function $w$ is entire if its domain is $\mathbb{C}$.
\begin{theorem}[ Liouville's Theorem, \cite{13}]\label{t2}
Let $w:\mathbb{C} \to \mathbb{C}$ be an bounded entire function. Then $w$ is constant.
\end{theorem}
Liouville's theorem states that every bounded entire function must be constant. That is, every holomorphic function $w$ for which there exists a positive number $M$ such that $|w(z)|\leq M$ for all $z$ in $\mathbb {C}$  is constant. Equivalently, non-constant holomorphic functions on $\mathbb {C}$  have unbounded images. The theorem is considerably improved by Picard's little theorem, which says that every entire function whose image omits two or more complex numbers must be constant. As a matter of fact, one can simplify Liouville's theorem to two equations given by
\[\left\{ \begin{matrix}
   \frac{\partial w}{\partial \overline{z}}=0  \\
   \left| w\left( z \right) \right|\le M,M>0  \\
\end{matrix} \right.\]
for $z$ in the whole domain $\mathbb {C}$.

The theorem follows from the fact that holomorphic functions are analytic. If $w$ is an entire function, it can be represented by its Taylor series about 0:
$w(z)=\sum _{k=0}^{\infty }a_{k}z^{k}$,
where
\[a_{k}={\frac {w^{(k)}(0)}{k!}}={1 \over 2\pi \sqrt{-1}}\oint _{C_{r}}{\frac {w(\zeta )}{\zeta ^{k+1}}}\,d\zeta\]
and $C_{r}$ is the circle about 0 of radius $r > 0$. Suppose $w$ is bounded: i.e. there exists a constant $M$ such that $|w(z)|\leq M$ for all $z$. We can estimate directly
\[|a_{k}|\leq {\frac {1}{2\pi }}\oint _{C_{r}}{\frac {|w(\zeta )|}{|\zeta |^{k+1}}}\,|d\zeta |\leq {\frac {M}{r^{k}}}\]
where in the second inequality we have used the fact that $|z|=r$ on the circle $C_{r}$. But the choice of $r$ in the above is an arbitrary positive number. Therefore, letting $r$ tend to infinity since $w$ is analytic on the entire plane gives $a_{k} = 0$ for all $k\geq 1$. Thus $w(z) = a_{0}$ and this proves the theorem.

If $w$ is less than or equal to a scalar times its input, then it is linear.
Suppose that $w$ is entire and $|w(z)|$ is less than or equal to $M|z|$, for
$M$ a positive real number. We can apply Cauchy's integral formula; we have that
\[|w'(z)|={\frac {1}{2\pi }}\left|\oint _{C_{r}}{\frac {w(\zeta )}{(\zeta -z)^{2}}}d\zeta \right|\leq {\frac {MI}{2\pi }}\]
where $I$ is the value of the remaining integral. This shows that $w' $ is bounded and entire, so it must be constant, by Liouville's theorem.

The following result corresponds to the definable analogue of Liouville's theorem proved by Peterzil and Starchenko \cite{13}.

\begin{theorem}[\cite{13}]
 Let $w:\mathbb{C} \to \mathbb{C}$ be a definable bounded entire function. Then $w$ is constant.
\end{theorem}
physical proof of Liouville’s theorem for a class generalized harmonic functions by the method of parabolic equation[see more \cite{12}].

\begin{lemma}[\cite{23}]\label{l5}
  If $w$ is an analytic function on $\Omega\in \mathbb{C}$ and if it satisfies one of conditions below,
  \begin{enumerate}
    \item $w'\left( z \right)=0$ holds on $\Omega$.
    \item $\left| w\left( z \right) \right|$ is a constant on $\Omega$.
    \item $\overline{w}\left( z \right)$ is analytic in $\Omega$.
    \item $\operatorname{Re}w\left( z \right)$, or $\operatorname{Im}w\left( z \right)$ is constant in $\Omega$..
  \end{enumerate}
  then $w$ is a constant on $\Omega$
\end{lemma}
Eventually, the condition 1 is the lemma \ref{l2}, here $K=1$ obviously.

\subsection{Maximum modulus principle}
\begin{theorem}[Maximum Modulus Principle, \cite{17,18}]
   Let $w:\Omega \to \mathbb{C}$, where $\Omega$ is open and
connected, be analytic. If $|w|$ has a local maximum, then $w$ is a constant.

\end{theorem}

the maximum modulus principle in complex analysis states that if $w$ is a holomorphic function, then the modulus $|w|$ cannot exhibit a true local maximum that is properly within the domain of $w$.

Let $w$ be a function holomorphic on some connected open subset $\Omega$ of the complex plane $\mathbb{C}$ and taking complex values. If $z_{0}$ is a point in $\Omega$ such that $|w(z_{0})|\geq |w(z)|$
for all $z$ in a neighborhood of $z_{0}$, then the function $w$ is constant on $\Omega$.  In other words, either $w$ is a constant function, or, for any point $z_{0}$ inside the domain of $w$ there exist other points arbitrarily close to $z_{0}$ at which $|w|$ takes larger values.  \cite{14} has given a sneaky proof of the maximum modulus principle.

The maximum modulus principle has many uses in complex analysis, and may be used to prove the following:
The fundamental theorem of algebra.
Schwarz's lemma, a result which in turn has many generalisations and applications in complex analysis.

\subsection{Cauchy's integral theorem}
the Cauchy integral theorem in complex analysis is an important statement about line integrals for holomorphic functions in the complex plane. Essentially, it says that if two different paths connect the same two points, and a function is holomorphic everywhere in between the two paths, then the two path integrals of the function will be the same \cite{19}.

The theorem is usually formulated for closed paths as follows:
\begin{theorem}[\cite{5}]
  let $U\subset\mathbb{C}$ be an open subset, let $w: U\rightarrow \mathbb{C}$ be a holomorphic function, and let $\gamma$  be a closed curve in an open set $U$. Then
\[\oint _{\gamma }w(z)\,dz=0\]
\end{theorem}
If one assumes that the partial derivatives of a holomorphic function are continuous, the Cauchy integral theorem can be proved as a direct consequence of Green's theorem and the fact that the real and imaginary parts of $w=u+\sqrt{-1}v $  must satisfy the Cauchy-Riemann equations in the region bounded by $\gamma$, and moreover in the open neighborhood $U$ of this region.

As was shown by \'{E}douard Goursat, Cauchy's integral theorem can be proven assuming only that the complex derivative $w′(z)$ exists everywhere in $U$. This is significant, because one can then prove Cauchy's integral formula for these functions, and from that deduce these functions are in fact infinitely differentiable.

Morera's theorem states that a continuous, complex-valued function $w$ defined on an open set $\Omega$ in the complex plane that satisfies
$\oint _{\gamma }w(z)\,dz=0$
for every closed piecewise $C_{1}$ curve $\gamma$  in $\Omega$ must be holomorphic.

The Cauchy integral theorem leads to Cauchy's integral formula and the residue theorem.

\begin{theorem}[\cite{5}]
Suppose $w$ is holomorphic in an open set that contains
the closure of a disc $\Omega$. If $C$ denotes the boundary circle of this disc with the positive orientation, then
\[w\left( z \right)=-\frac{\sqrt{-1}}{2\pi }\int\limits_{C}{\frac{w\left( \zeta  \right)}{\zeta -z}}d\zeta \] for any point $z\in \Omega$.
\end{theorem}
The proof of this statement uses the Cauchy integral theorem and like that theorem it only requires $w$ to be complex differentiable.

Cauchy's integral formula is a central statement in complex analysis. It expresses the fact that a holomorphic function defined on a disk is completely determined by its values on the boundary of the disk, and it provides integral formulas for all derivatives of a holomorphic function. Cauchy's formula shows that, in complex analysis, differentiation is equivalent to integration: complex differentiation, like integration, behaves well under uniform limits a result denied in real analysis.

The theorem stated above can be generalized. The circle $\gamma$ can be replaced by any closed rectifiable curve in $U$ which has winding number one about $z$. Moreover, as for the Cauchy integral theorem, it is sufficient to require that $w$ be holomorphic in the open region enclosed by the path and continuous on its closure.

Note that not every continuous function on the boundary can be used to produce a function inside the boundary that fits the given boundary function.

\begin{lemma}[Cauchy's Estimate, \cite{5,6}]\label{l4}
  Let $w$ be analytic in $B(a; R)$ and suppose if $\left| w\left( z \right) \right|\le M$  for all $z$ in $B(a; R)$. Then
  \[\left| {{w}^{\left( n \right)}}\left( a \right) \right|\le \frac{n!M}{{{R}^{n}}}\]
\end{lemma}

\begin{theorem}[\cite{5,6}]
Let $U\subset \mathbb{C}$ be an open set and let $w$ be holomorphic on $U$. Then $w\in {{C}^{\infty }}\left( U \right)$. Moreover, if $\overline{\Omega}\left( p,r \right)\subset U$ and $z\in \Omega\left( p,r \right)$, then
  \[\frac{{{\partial }^{k}}}{\partial {{z}^{k}}}w\left( z \right)=-\frac{k!\sqrt{-1}}{2\pi }\oint\limits_{\left| \zeta -p \right|=r}{\frac{w\left( \zeta  \right)}{{{\left( \zeta -z \right)}^{k+1}}}d\zeta },~~k=0,1,2,\cdots \]

\end{theorem}
This formula is sometimes referred to as Cauchy's differentiation formula.
The integral formula has broad applications. First, it implies that a function which is holomorphic in an open set is in fact infinitely differentiable there. Furthermore, it is an analytic function, meaning that it can be represented as a power series. The proof of this uses the dominated convergence theorem and the geometric series applied to
\[w(\zeta) = \frac{1}{2\pi \sqrt{-1}}\int_C \frac{w(z)}{z-\zeta}\,dz\]
The formula is also used to prove the residue theorem, which is a result for meromorphic functions, and a related result, the argument principle.  In addition the Cauchy formulas for the higher order derivatives show that all these derivatives also converge uniformly.

A version of Cauchy's integral formula is the Cauchy-Pompeiu formula, and holds for smooth functions as well, as it is based on Stokes' theorem. Let $\Omega$ be a disc in $\mathbb{C}$ and suppose that $w$ is a complex-valued $C^{1}$ function on the closure of $\Omega$. Then
\[{\displaystyle w(\zeta )={\frac {1}{2\pi \sqrt{-1}}}\int _{\partial \Omega}{\frac {w(z)dz}{z-\zeta }}-{\frac {1}{\pi }}\iint _{\Omega}{\frac {\partial w}{\partial {\bar {z}}}}(z){\frac {dx\wedge dy}{z-\zeta }}.}\]
One may use this representation formula to solve the inhomogeneous Cauchy-Riemann equations in $\Omega$. Indeed, if $\varphi$ is a function in $\Omega$, then a particular solution $w$ of the equation is a holomorphic function outside the support of $\mu$. Moreover, if in an open set $\Omega$,
$d\mu = \frac{1}{2\pi \sqrt{-1}}\phi \, dz\wedge d\bar{z}$
for some $\varphi \in C^{k}(\Omega), (k\geq 1)$, then $w(\zeta,\bar{\zeta})$ is also in $C^{k}(\Omega)$ and satisfies the equation
$ \frac{\partial w}{\partial\bar{z}} = \phi(z,\bar{z})$.
Note that for smooth complex-valued functions $w$ of compact support on $\mathbb{C}$ the generalized Cauchy integral formula simplifies to
\[w(\zeta) =  \frac{1}{2\pi \sqrt{-1}}\iint \frac{\partial w}{\partial \bar{z}}\frac{dz\wedge d\bar{z}}{z-\zeta}\]

\section{Generalized Cauchy's integral theorem}
Next we study the integral of the structural analytic function in the region $\mathbb{C}$. We know that there is a Green formula in ${{\mathbb{R}}^{2}}$  showing
\[\int\limits_{\partial \omega }{Pdx+Qdy}=\int\limits_{\omega }{\left( {{Q}_{x}}-{{P}_{y}} \right)d\lambda }\]
where $\omega$  is a bounded region whose boundary is composed of one or several smooth curves and taking right underneath.
If complex coordinate is used, the Green formula can be written as\cite{11}
\[\int\limits_{\partial \omega }{fdz+gd\overline{z}}=2\sqrt{-1}\int\limits_{\omega }{\left( \frac{\partial f}{\partial \overline{z}}-\frac{\partial g}{\partial z} \right)d\lambda \left( z \right)},~~f,g\in {{C}^{1}}\left( \overline{\omega } \right)\]Actually, for $f,g\in {{C}^{1}}\left( \overline{\omega } \right)$, setting $f,g$ are real functions, Otherwise, the real part and the imaginary part can be considered respectively. Let
$P=f+g,Q=\sqrt{-1}\left( f-g \right)$, plugging it into Green's formula, then
\[\int\limits_{\partial \omega }{fdz+gd\overline{z}}=2\sqrt{-1}\int\limits_{\omega }{\left( \frac{\partial f}{\partial \overline{z}}-\frac{\partial g}{\partial z} \right)d\lambda \left( z \right)}\]
In particular, if $g=0$ is taken, then
\[\int\limits_{\partial \omega }{fdz}=2\sqrt{-1}\int\limits_{\omega }{\frac{\partial f}{\partial \overline{z}}d\lambda \left( z \right)}\]
\begin{theorem}[Cauchy's theorem]
  If $f\in Hol\left( \omega  \right)\cap {{C}^{1}}\left( \overline{\omega } \right)$, then
  \[\int\limits_{\partial \omega }{fdz}=0\]

\end{theorem}

In this section, one will show that structural holomorphic function also can get Cauchy's integral theorem satisfied, it's now called generalized Cauchy's integral theorem. In the beginning, one knows that
\[\widetilde{w}=Kw=\left( {{k}_{1}}u-{{k}_{2}}v \right)+\sqrt{-1}\left( {{k}_{2}}u+{{k}_{1}}v \right)\]
together with $dz=dx+\sqrt{-1}dy$ can deduce the result as follows
\[\widetilde{w}dz=\left( {{k}_{1}}u-{{k}_{2}}v \right)dx-\left( {{k}_{2}}u+{{k}_{1}}v \right)dy+\sqrt{-1}\left[ \left( {{k}_{2}}u+{{k}_{1}}v \right)dx+\left( {{k}_{1}}u-{{k}_{2}}v \right)dy \right]\]

\begin{theorem}
 let $U\subset\mathbb{C}$ be an open subset, let $w: U\rightarrow \mathbb{C}$ be a structural holomorphic function, and let $\gamma$  be a closed curve in an open set $U$. Then
\[\oint\limits_{\gamma }{\widetilde{w}\left( z \right)dz}=0\]
where $\widetilde{w}=Kw$.
  \begin{proof}

In this case we have
\[\oint\limits_{\gamma }{\widetilde{w}dz}=2\sqrt{-1}\int\limits_{\omega }{\frac{\partial \widetilde{w}}{\partial \overline{z}}d\lambda \left( z \right)}=2\sqrt{-1}\int\limits_{\omega }{\frac{{\rm{D}}w}{\partial \overline{z}}d\lambda \left( z \right)}\]where $d\lambda \left( z \right)$ is area element.
Substituting the theorem \ref{t3} in 1-dimensional into above equation, and then it gives the desired result.

  \end{proof}

\end{theorem}

\begin{theorem}
  Suppose $w\in Shol\left( \Omega \right)$ is structural holomorphic in an open set that contains
the closure of a disc $\Omega$. If $C$ denotes the boundary circle of this disc with
the positive orientation, then
\[\widetilde{w}\left( z \right)=-\frac{\sqrt{-1}}{2\pi }\int\limits_{C}{\frac{\widetilde{w}\left( \zeta  \right)}{\zeta -z}}d\zeta \]
for any point $z\in \Omega$.

\end{theorem}

\section{Liouville's theorem on structural holomorphic}
\begin{lemma}\label{l3}
Suppose that $w\in Shol\left( \Omega \right)$ is a structural holomorphic function on a connected open set
$\Omega \subset \mathbb{C}$. If $\frac{\rm{D}}{\partial z}w=0$ on $\Omega$, then $w=\Phi {{e}^{-K}}$, where $\Phi =Const$.
  \begin{proof}
    Since $w$ is structural holomorphic by the theorem, then $\frac{{\rm{D}}}{\partial \overline{z}}w=0$. Therefore, we only need to solve the system of equation below
\begin{align}
  & \frac{\rm{D}}{\partial z}w=0 \notag\\
 & \frac{\rm{D}}{\partial \overline{z}}w=0 \notag
\end{align}
Equivalently, one has to get solution from equation
\[\frac{\partial w}{\partial \overline{z}}=-w\frac{\partial K}{\partial \overline{z}},~~~\frac{\partial w}{\partial z}=-w\frac{\partial K}{\partial z}\]
Its solution can be expressed as \[w=\Phi {{e}^{-K}}\]
where $\frac{\partial }{\partial z}\Phi =0,~\frac{\partial }{\partial \overline{z}}\Phi =0$. According to the lemma \ref{l2} and lemma \ref{l5}, it's easy to get a result that $\Phi =Const$.

  \end{proof}

\end{lemma}
In fact, for the formal solution, we get \[w=\Phi {{e}^{-K}}=\Phi \left( 1-K+\cdots  \right)=\Phi +U\left( K \right)\]
where $U\left( K \right)$ is a complex function in terms of the $K$, obviously, $U\left( 0\right)=0$ holds for the classical result.  Conversely, it gets ${{e}^{K}}w=\Phi$, and \[\left| {{e}^{K}}w \right|=\left| {{e}^{K}} \right|\left| w \right|={{e}^{{{k}_{1}}}}\left| w \right|=\left| \Phi  \right|\]thusly, it obtains \[\left| w \right|\left( {{k}_{1}} \right)=\left| \Phi  \right|{{e}^{-{{k}_{1}}}}\]It clearly implies that the module of the complex function $w$ is the function of the real function ${k}_{1}$, it reveals that there exists maximum value of the module of $w$ when ${k}_{1}=0$ is taken such that ${{\left| w \right|}_{Max}}\left( 0 \right)=\left| \Phi  \right|$. Subsequently, it leads to the restriction
\begin{equation}\label{eq2}
  \left| w \right|\left( {{k}_{1}} \right)=\left| \Phi  \right|{{e}^{-{{k}_{1}}}}\le {{\left| w \right|}_{Max}}\left( 0 \right)=\left| \Phi  \right|
\end{equation}
Accordingly \[\left| w \right|\left( {{k}_{1}} \right)={{\left| w \right|}_{Max}}\left( 0 \right){{e}^{-{{k}_{1}}}}\]
As a consequence, it implies that $w$ is always a bounded function in  $\mathbb{C}$.

Actually, as a complex function is a structural entire function,  it means
  $\frac{{\rm{D}}}{\partial \overline{z}}w=0$ in $\mathbb{C}$, then it gets a formal solution shown as $w=\Phi {{e}^{-K}}$ in which $\frac{\partial }{\partial \overline{z}}\Phi =0$ means that $\Phi$ is an arbitrary entire
function. Accordingly, it automatically satisfies the bound \eqref{eq2}.   Thusly, it remains the consistency to the structural Liouville's theorem below.

\begin{corollary}[Structural Liouville's theorem]\label{c2}
Let $w:\mathbb{C} \to \mathbb{C}$ be an bounded structural entire function. Then $w=\Phi {{e}^{-K}}$, where $\Phi =Const$.
\end{corollary}
Note that corollary \ref{c2} is degenerated to the theorem \ref{t2} if $K=Const$ holds for $z$ in the whole domain $\mathbb {C}$.

\subsection{Maximum modulus principle on structural holomorphic }

Let $w$ be a function structural holomorphic on some connected open subset $\Omega$ of the complex plane $\mathbb{C}$. Based on the lemma \ref{l3} and corollary \ref{c2}, we can see that as indicated, if $w(z)$ is a structural holomorphic entire function, then it gets a formal solution shown as $w=\Phi {{e}^{-K}}$ with $\frac{\partial }{\partial \overline{z}}\Phi =0$, its module is $\left| w \right|=\left| \Phi  \right|{{e}^{-{{k}_{1}}}}$ that is a function in terms of ${k}_{1}$, clearly, its maximum is obtained only if ${k}_{1}=0$, then $\left| w \right|=\left| \Phi  \right|$, it implies that the maximum of $\left| w \right|$ under the condition of structural holomorphic connects to the function ${k}_{1}$ and the module of $\Phi$, notice that the function ${k}_{1}$ is only associated with the complex domain, it directly has impact on the maximum of $\left| w \right|$, however, in general case, ${k}_{1}\neq 0$ always holds.

\begin{corollary}
  Let $w(z)$ be a structural holomorphic entire function, and $|w|$ takes larger values at $z_{0}$. then $\left| w \right|({k}_{1})=\left| \Phi  \right|{{e}^{-{{k}_{1}}}}$ in $\mathbb{C}$, and ${{\left| w \right|}_{Max}}\left( 0 \right)=\left| \Phi  \right|$ for the classical maximum modulus principle.
\end{corollary}
Note that when $K$ is a constant function, the maximum modulus principle still holds in $\mathbb{C}$, but as $K$ is a function in terms of the variables $z$ or $x,y$, the maximum modulus principle will be accordingly rewritten in a way of structural function $K$. As it shown above, it relies on the magnitude of structural function $K$. It is well-known that the maximum modulus principle in complex analysis holds at $K=Const$ only.

\end{document}